\documentclass{jT}
\usepackage{amsmath}
\usepackage{amssymb}
\usepackage{amsfonts}
\usepackage{latexsym}

\issueinfo{15}% volume number
  {2003}%     % year
\makeatletter
\def\start@page{627}
\def\end@page{637}
\makeatother

\date{3 d\'ecembre 2001}

\def\cal{\mathcal}
\def\iso{\cong}
\def\gl{\mathop{\rm GL}\nolimits}
\def\mult{\mathop{\rm mult}\nolimits}
\def\cross{\times}
\def\inv{^{-1}}
\def\frob{{\mathop{\rm Fr}\nolimits}}
\def\cala{{\cal A}}

\def\calc{{\cal C}}

\def\calf{{\cal F}}
\def\calg{{\cal G}}

\def\calm{{\cal M}}

\def\calp{{\cal P}}

\def\jac{\mathop{\rm Jac}\nolimits}

\def\aut{\mathop{\rm Aut}\nolimits}
\def\spec{\mathop{\rm Spec}\nolimits}
\def\sriso{\stackrel{\sim}{\rightarrow}}

\def\diag{\mathop{\rm diag}\nolimits}
\newcommand{\half}[1]{\frac{#1}{2}}
\def\ff{{\mathord{\mathbb{F}}}}
\def\gp{{\mathord{\mathbb{G}}}}
\def\cx{{\mathord{\mathbb{C}}}}

\def\aff{{\mathord{\mathbb{A}}}}

\def\ra{\rightarrow}
\newcommand{\sra}[1]{\stackrel{#1}{\ra}}
\def\ndiv{\nmid}
\def\gal{\mathop{\rm Gal}}
\newcommand{\oneover}[1]{\frac{1}{#1}}
\newcommand{\twomat}[4]{\begin{array}{ll} {#1} & {#2} \\ {#3} & {#4}
\end{array}}
\def\id{\mathop{\rm id}\nolimits}
\newcommand{\abs}[1]{{\left|#1\right|}}
\newcommand{\rest}[1]{|_{#1}}

\newcommand{\st}[1]{\{#1\}}
\newcommand{\ip}[1]{{{\langle #1 \rangle}}}

\newcommand{\mt}[1]{\left(#1\right)}

\newcommand{\twmt}[4]{\mt{\twomat{#1}{#2}{#3}{#4}}}
\def\units{^\cross}

\def\sp{{\mathop{\rm Sp}}}

\def\gsp{{\mathop{\rm GSp}}}

\def\geom{{\rm geom}}

\def\eigenone{{\sf (E)}}
\def\notone{{\sf (N)}}
\def\eigenweird{{\sf (R)}}
\def\fm{{\sf(FM)}}

\input{diagrams}
\newarrow{to}---->
\newarrow{mto}|--->
\newarrow{impto}===={=>}
\newarrowtail{<=}\Rightarrow\Leftarrow\Uparrow\Downarrow
\newarrow{bboth}{<=}==={=>}
\newarrow{inject}{hooka}---{vee}
\newarrow{surject}----{>>}

\newtheorem{theorem}{Theorem}[section]
\newtheorem{lemma}[theorem]{Lemma}

\newtheorem{prop}[theorem]{Proposition}

\newtheorem{cor}[theorem]{Corollary}

\newtheorem{conj}[theorem]{Conjecture}
\newtheorem{question}[theorem]{Question}

\theoremstyle{definition}
\newtheorem*{definition}{Definition}
\newtheorem*{remark}{Remark}

\begin{document}
\title[Notes on an analogue of the Fontaine-Mazur conjecture]{Notes on an analogue of the Fontaine-Mazur conjecture}

\author[Jeffrey D. {\sc Achter}]{{\sc Jeffrey D.} ACHTER}
\address{Jeffrey D. {\sc Achter}\\
Department of Mathematics\\
Columbia University\\
New York, NY 10027, USA}
\email{achter@math.columbia.edu}

\author[Joshua {\sc Holden}]{{\sc Joshua} HOLDEN}
\address{Joshua {\sc Holden}\\
Department of Mathematics\\
Rose-Hulman Institute of Technology\\
Terre Haute, IN 47803, USA}
\email{holden@rose-hulman.edu}

\maketitle

\begin{resume}
On estime le proportion des corps de fonctions qui remplissent des
conditions qui impliquent un analogue de la conjecture de Fontaine et Mazur.
En passant, on calcule le proportion des vari\'et\'es abeli\'ennes (ou
Jacobiennes) sur un corps fini qui poss\`edent un point rationnel d'orde
$\ell$.
\end{resume}

\begin{abstr}
   We estimate the proportion of function fields satisfying certain
   conditions which imply a function field analogue of the
   Fontaine-Mazur conjecture.  As a byproduct, we compute the fraction
   of abelian varieties (or even Jacobians) over a finite field which
   have a rational point of order $\ell$.
\end{abstr}

\bigskip

\section{Introduction}

The paper~\cite{Holden99} discusses the following conjecture,
originally stated by Fontaine and Mazur in~\cite{FM}:

\begin{conj}[Fontaine-Mazur, as restated in~\cite{Boston98}]
Let $F$ be a number field and $\ell$ any prime.
There does not exist an infinite everywhere unramified Galois pro-$\ell$
extension $M$ of $F$ such that
$\gal(M/F)$ is uniform.
\end{conj}

The definitions of \emph{powerful} and \emph{uniform} are
taken from~\cite{DDMS}:

\begin{definition}
Let $G$ be a pro-$\ell$ group.  $G$ is \emph{powerful}
if $\ell$ is odd and $G/\overline{G^{\ell}}$ is abelian, or if
$\ell = 2$ and $G/\overline{G^{4}}$ is abelian.  ($G^{n}$ is the
subgroup of $G$ generated by the $n$-th powers of elements in $G$,
and $\overline{G^{n}}$ is its closure.)
\end{definition}

\begin{definition}
A pro-$\ell$ group is \emph{uniformly powerful}, or just \emph{uniform},
if (i)~$G$ is finitely generated, (ii)~$G$ is powerful, and (iii)~for
all $i$, $\left[\overline{G^{\ell^{i}}} : \overline{G^{\ell^{i+1}}}
\right] = \left[G : \overline{G^{\ell}} \right]$.
\end{definition}

The paper~\cite{Holden99} then raises the following question:
\begin{question}
Let $F$ be a function field over a finite field $k_0$, and $\ell$
 a prime invertible in $k_0$.  Let  $F_\infty =
Fk_{\ell^\infty}$ be obtained  by taking the maximal
pro-$\ell$ extension of the constant field.  Does $F$ satisfy the
following property?

\begin{center}
\fm \ \ \
\parbox{3in}{
Let $F'$ be any non-trivial subextension of $F_\infty/F$, and
$M$ be any infinite unramified pro-$\ell$
 extension of $F'$.  If $M$ is  Galois over $F$ and $M$ does not contain
$F_\infty$, then $\gal(M/F')$ cannot be uniform.}
\end{center}
\end{question}

(See~\cite{Boston98} for a discussion of the relationship between the
conjecture as phrased here and the conjecture as originally given.)

The general answer to this question is in fact negative, as shown by
examples due to Ihara \cite{Ihara83} and to Frey, Kani, and
V\"olklein \cite{FKV}.  In fact, there is reason to believe that the
correct analogue of the Fontaine-Mazur Conjecture will be found not in
questions related to~\cite{Holden99} but in work related to that of
de~Jong~\cite{deJong}.

Nevertheless, results of \cite{Holden99} answer the above question
affirmatively in a large class of situations (see \ref{oldmainth}
below).  Since a great deal of effort has been put into constructing
fields which do not satisfy \fm, we would like to know if they are in
fact common, or if they are rather rare.  The present paper will
attempt to quantify in some way the proportion of fields $F$ which
satisfy \fm.

The strategy is simple enough.  The aforementioned paper
\cite{Holden99} provides conditions on $F$ (such as the absence of
an $\ell$-torsion  element in the class group of $F$; see
\ref{oldmainth})
which force an affirmative
answer to the question.  These conditions may be formulated in terms
of the action of Frobenius on the $\ell$-torsion of the Jacobian of
the smooth, proper model of $F$.  Equidistribution results for
$\ell$-adic monodromy imply analogous results for $\bmod\ \ell$
monodromy, and show that Frobenius automorphisms are evenly distributed among
$\gsp_{2g}(\ff_\ell)$; counting symplectic similitudes then finishes
the analysis.  As a pleasant side effect, we calculate the proportion
of abelian varieties over $k_0$ with a $k_0$-rational point of order
$\ell$.

Section \ref{mono} reviews work of Katz on equidistribution, and
axiomatizes our situation.  Section \ref{sp} studies (the size of)
certain conjugacy classes in $\gsp_{2g}(\ff_\ell)$.  The final section
gives the quantitative Fontaine-Mazur results alluded to in the title.

This paper was written while the first author visited the second at
Duke University as part of the second DMJ-IMRN conference; we thank
these institutions for providing such a pleasant working environment.
We also thank our anonymous referee for helpful suggestions.

\section{Monodromy groups}\label{mono}
The main piece of technology which drives this paper is an
equidistribution theorem for lisse $\ell$-adic sheaves.  Originally due
to Deligne \cite[3.5]{delw2}, it has since been clarified and
amplified by Katz.  Deferring to chapter nine of \cite{katzsarnak} a
careful and complete exposition of the theory,  we content ourselves
by recalling the precise result needed here.

Let $(V,\ip{\cdot,\cdot})$ be a $2g$-dimensional vector space over
$\ff_\ell$ equipped with a symplectic form.  Recall the definition of
the group of symplectic similitudes of $(V,\ip{\cdot,\cdot})$:
\begin{multline*}
\gsp_{2g}(\ff_\ell) \\
\begin{split}
&\iso \gsp(V,\ip{\cdot,\cdot}) \\
&= \st{
A \in \gl(V) |  \exists \mult(A)\in \ff_\ell\units\text : \forall v,w\in V,
\ip{Av,Aw} = \mult(A)\ip{v,w}}.
\end{split}
\end{multline*}
The ``multiplicator'' $\mult$ is a character of the group, and its
kernel is the usual symplectic group $\sp_{2g}(\ff_\ell)$.  For
$\gamma\in \ff_\ell\units$, let $\gsp_{2g}^\gamma(\ff_\ell) =
\mult\inv(\gamma)$ be the
set of symplectic similitudes with multiplier $\gamma$.  Each
$\gsp_{2g}^\gamma$ is a torsor over $\sp_{2g}$.

Now let $k_0 = \ff_q$ be a finite field of characteristic $p$, prime to
$\ell$, and let $U/k_0$ be a smooth, geometrically irreducible variety
with geometric generic point $\bar\eta$.
If $k$ is a finite extension of $k_0$, then one may associate to any
point $u\in U(k)$ its (conjugacy class of) Frobenius $\frob_{u/k}$ in
$\pi_1(U,\bar\eta)$.

%Suppose $\calf$ is a local system of symplectic $V$-modules on $U$.
%Such an object is tantamount to a representation $\rho_\calf: \pi_1(U,\bar\eta)
%\ra \aut(\calf_{\bar\eta}) \iso \gsp_{2g}(\ff_\ell)$.

Suppose $\calf$ is a local system of symplectic $\ff_\ell$-modules of
rank $2g$ on $U$.  Recall that such an object is tantamount to a
continuous
representation $\rho_\calf: \pi_1(U,\bar\eta)
\ra \aut(\calf_{\bar\eta}) \iso \gsp_{2g}(\ff_\ell)$.  (To see this,
one may consider the total space of $\calf$, which is an
\'etale cover of $U$.  The fundamental group of $U$ acts on covering
spaces of $U$, and in particular on the total space of $\calf$; this
is the desired representation.)

A simple case of Katz's equidistribution theorem says the following:

\begin{theorem}[Katz]
\label{katz}
In the situation above, suppose the sheaf gives rise to a commutative
diagram
\begin{diagram}
1 & \rto & \pi_1^\geom(U,\bar\eta) & \rto  & \pi_1(U,\bar\eta) & \rto & \gal(\bar{k_0}/k_0) & \rto
& 1 \\
&&\dsurject>{\rho^\geom} &&\dto>\rho &&\dto>{\rho^{k_0}}  \\
1 & \rto & \sp_{2g}(\ff_\ell) & \rto & \gsp_{2g}(\ff_\ell) & \rto & \gp_m(\ff_\ell) & \rto & 1
\end{diagram}
where $\rho^\geom$ is surjective.  There
is a constant $C$ such that, for any union of conjugacy classes
$W\subset \gsp_{2g}(\ff_\ell)$ and any finite extension $k$ of $k_0$,
\[
\abs{
\frac{
\#\st{ u \in U(k) : \rho(\frob_{u,k}) \in W } }
{ \#U(k)}
-
\frac{ \#(W\cap \gsp_{2g}^{\gamma(k)}(\ff_\ell))}{\sp_{2g}(\ff_\ell)}}
\le
\frac C{\sqrt{\#k}},
\]
where $\gamma(k)$ is the image of the canonical generator of
$\gal(k)$ under $\rho^{k_0}$.
\end{theorem}

\begin{proof}
This is a special case of \cite[9.7.13]{katzsarnak}; see also
\cite[4.1]{chav}.
\end{proof}

The constant $C$ in Katz's theorem is effectively computable, and the
theorem actually holds uniformly in families; but we will not need
such developments here.
\def\cgn{{_N\calc_g}}
\def\mgn{{_N\calm_g}}
\def\cgln{{_{\ell N}\calc_g}}
\def\mgln{{_{\ell N}\calm_g}}

Let $\calc\ra\calm \ra \spec k_0$ be a smooth, irreducible family of
curves of genus $g \ge 1$.  There is a sheaf $\calf = \calf_{\calc,\ell}$ of
abelian groups on $\calm$ whose fiber at a geometric point $\bar x\in
\calm$ is the $\ell$-torsion of the Jacobian $\jac(\calc_x)[\ell]$.
We will say that the family of curves has full $\ell$-monodromy if the
associated representation $\rho_\calf: \pi_1^\geom(\calm, \bar\eta)
\ra \sp_{2g}(\ff_\ell)$ is surjective.  In practice, general families
of curves tend to have full $\ell$-monodromy; see, for instance, the
introduction to \cite{ekedahl}.  Concretely, we will see below that
the universal family of curves over $\mgn$ has full $\ell$-monodromy.

\begin{lemma}\label{full} Let $C \ra S \ra \spec k_0$ be a geometrically
irreducible, smooth versal family of proper smooth curves of genus
$g$.  For almost all $\ell$, $C/S$ has full $\ell$-monodromy.
\end{lemma}

\begin{proof}
Fix a natural number $N$ relatively prime to $p$ and consider $\cgn\ra
\mgn$, the universal curve of genus $g$ with principal Jacobi
structure of level $N$.  If $\ell | N$, then the final paragraphs of
\cite{dm} imply that this family has full $\ell$-monodromy.  Indeed,
\cite[5.11]{dm} shows that it suffices to verify the statement for the
analogous family over $\cx$, and
\cite[5.13,5.15]{dm} provides this proof.  If $\ell$ is relatively
prime to $N$, then consider the moduli space $\mgln$.  On one hand,
the $\ell$-torsion of the Jacobian of $\cgln$ has full
$\ell$-monodromy.  On the other hand, the forgetful map $\mgln \ra
\mgn$ is finite; therefore, $\cgn\ra \mgn$ has full $\ell$-monodromy,
too.

For any $C/S$ as in the statement of the lemma,  there is an \'etale base change
$T\sra\phi S$ so that $\phi^*T$ admits a level $N$ structure.
Then $\phi^*T$ is the pullback of $\cgn$ by the classifying map $\psi:
T \ra \mgn$.  Moreover, the sheaf of $\ell$-torsion on $\phi^*T$,
$\calf_{\phi^*T,\ell}$ is the pullback of the universal $\ell$-torsion:
\begin{diagram}
\calf_{\phi^*T,\ell} & \rto & \calf_{\cgn,\ell} \\
\dto&&\dto \\
T & \rto^\psi & \mgn
\end{diagram}

By the versality assumption, $T$ has dense image in $\mgn$.
We have seen above that $\calf_{\cgn,\ell} \ra \mgn$ has monodromy
group $\sp_{2g}(\ff_\ell)$.  Thus, as long as $\ell\ndiv \deg
\psi$, $\calf_{\phi^*T,\ell} \ra T$ has full $\ell$-monodromy; {\em a
fortiori}, $C/S$ does, too.
\end{proof}

We now relate these notions to the quantitative Fontaine-Mazur
question posed at the beginning of this paper.  Let $\calp$ be a
property of abelian varieties over finite extensions $k$ of $k_0$
which is detectable on $\ell$-torsion.  We will say a curve has
$\calp$ if its Jacobian does, and that a function field has $\calp$ if
its smooth, proper model does.  (We have in mind, e.g., having a
$k$-rational point of order $\ell$.)  Define $W_\calp\subset
\gsp_{2g}(\ff_\ell)$ as the set of Frobenius automorphisms satisfying
$\calp$, and let $W_\calp^\gamma = W_\calp \cap \gsp_{2g}^\gamma$.
More precisely, $w\in W_\calp$ if and only if there exists an abelian
variety $X/k$ over a finite extension of $k_0$ and an isomorphism
$(V,w) \sriso (X[l](\bar k), \frob_{X/k})$.

\begin{cor} Fix a $\gamma\in \ff_\ell\units$.  Let $\st{k_n}$ be a collection
of extensions of $k_0$ such that $\lim_{n\ra\infty} \#k_n = \infty$,
and, for all $n$, $\gamma(k_n) = \gamma$.  Suppose that $\calc \ra
\calm \ra k_0$ is a smooth, irreducible family of curves with full
$\ell$-monodromy.  If $\calp$ is a property as
above, then
\[
\lim_{n\ra\infty}
\frac{\#\st{x\in \calm(k) : \calp(\calc_x)}}
{\#\calm(k)}
=
\frac{\# W^\gamma_\calp(\ff_\ell)}
{\#\sp_{2g}(\ff_\ell)}.
\]
\end{cor}

\begin{proof}
In view of the preceding discussion, this is an immediate
application of \ref{katz}.
\end{proof}

Let $\Xi^\gamma_g$ denote the set of all characteristic polynomials of
elements of $\gsp_{2g}^\gamma$.  It is well-known that $\Xi^\gamma_g
\iso \aff^g$; the isomorphism is given by sending a characteristic
polynomial to its first $g$ coefficients.  For a property $\calp$ as
above, let $\Psi^\gamma_\calp$ denote the set of all characteristic
polynomials which satisfy $\calp$.  The proportion of characteristic
polynomials satisfying $\calp$ is roughly the same as the proportion
of actual elements of $\gsp_{2g}^\gamma$ satisfying $\calp$.

\begin{lemma}
\label{psitow} For any property as above,
\[
\left(
\frac{\ell}{\ell+1}
\right)^{2g^2+g}
\frac{\#\Psi_{\calp,g}^\gamma(\ff_\ell)}{\#\Xi_g^\gamma(\ff_\ell)}
\le
\frac{\#W_\calp^\gamma(\ff_\ell)}{\#\sp_{2g}(\ff_\ell)}
\le
\left(
\frac{\ell}{\ell-1}
\right)^{2g^2+g}
\frac{\#\Psi_{\calp,g}^\gamma(\ff_\ell)}{\#\Xi_g^\gamma(\ff_\ell)}.
\]
\end{lemma}

\begin{proof}
For $f(x) \in \Xi^\gamma_g$, let $\Delta(f)$ be the number of
elements of $\gsp_{2g}^\gamma(\ff_\ell)$ whose characteristic
polynomial is $f(x)$.  One knows \cite[3.5]{chav} that, since $\dim
\sp_{2g} = 2g^2+g$,
\[
\frac{\ell^{2g^2}\#\sp_{2g}(\ff_\ell)}
{(\ell+1)^{2g^2+g}}
\le
\Delta(f)
\le
\frac{\ell^{2g^2}\#\sp_{2g}(\ff_\ell)}
{(\ell-1)^{2g^2+g}}.
\]
Adding up over all
elements of $W^\gamma_g$ we see that $\#W^\gamma_{\calp,g}(\ff_\ell)
 = \sum_{f\in \Xi^\gamma_{\calp,g}} \Delta(f)$, and thus
\begin{eqnarray*}
\frac{\ell^{2g^2}\#\sp_{2g}(\ff_\ell)}
{(\ell+1)^{2g^2+g}}\#\Psi^\gamma_{\calp,g}(\ff_\ell)
&\le
\#W^\gamma_{\calp,g}(\ff_\ell)
\le
&
\frac{\ell^{2g^2}\#\sp_{2g}(\ff_\ell)}
{(\ell-1)^{2g^2+g}} \#\Psi^\gamma_{\calp,g}(\ff_\ell); \\
\frac{\ell^{2g^2}}
{(\ell+1)^{2g^2+g}}\#\Psi^\gamma_{\calp,g}(\ff_\ell)
&\le
\frac{\#W^\gamma_{\calp,g}(\ff_\ell)}{\#\sp_{2g}(\ff_\ell)}
\le
&
\frac{\ell^{2g^2}}
{(\ell-1)^{2g^2+g}} \#\Psi^\gamma_{\calp,g}(\ff_\ell); \\
\left(\frac{\ell}
{\ell+1}\right)^{2g^2+g}
\frac{\#\Psi^\gamma_{\calp,g}(\ff_\ell)}{\#\Xi^\gamma_g(\ff_\ell)}
&\le
\frac{\#W^\gamma_{\calp,g}(\ff_\ell)}{\#\sp_{2g}(\ff_\ell)}
\le
&
\left(\frac{\ell}
{\ell-1}\right)^{2g^2+g}
\frac{\#\Psi^\gamma_{\calp,g}(\ff_\ell)}{\#\Xi^\gamma_g(\ff_\ell)}.
\end{eqnarray*}
\end{proof}

\section{Remarks on symplectic groups}\label{sp}

\subsection{Eigenvalue one}

We start by counting the number of
matrices for which $1$ is an eigenvalue; these will correspond to a
certain class of function fields which we will later show (in Theorem~\ref{oldmainth})
satisfy \fm.
Let $\eigenone$ be the
property of having $1$ as an eigenvalue.  Writing $f_A(x)$ for the
characteristic polynomial of $A\in \gsp_{2g}(\ff_\ell)$, we see that
$A\in W_{\eigenone,g}$ if and
only if $f_A(1) = 0$.   Barring any obvious reason to the contrary, one
might suppose that the values $\st{f_A(1)}_{A\in \gsp_{2g}(\ff_\ell)}$
are evenly
distributed in $\ff_\ell$, and thus that
$\#W_{\eigenone,g}/\#\gsp_{2g}(\ff_\ell)$ is about
$\oneover\ell$.  We will now show that this rough estimate is the
approximate truth -- and that, confounding our initial expectations,
$\oneover{\ell-1}$ is an even better estimate.

We need a little more notation in order to state our result precisely.  Let
\mbox{$T(g,\gamma,\ff_\ell) = \#W_{\eigenone, g}^\gamma$}  be the number of
elements of
$\gsp^\gamma_{2g}(\ff_\ell)$ which have one as an eigenvalue.  If $\gamma
\not = 1$, let $S(g,\gamma,\ff_\ell)$ be the number of
elements of $\gsp^\gamma_{2g}(\ff_\ell)$ for which the eigenspace
associated to $1$ has dimension $g$, and let $S(g,1,\ff_\ell)$ be the number of
unipotent symplectic matrices of rank $2g$.  Our goal in this section
is to compute $T(g,\gamma,\ff_\ell)/\#\gsp^\gamma_{2g}(\ff_\ell)$.  As an
organizational tool, we collect intermediate results in a series of
easy lemmas.

\begin{lemma} For $\gamma\in \ff_\ell\units$ and $S$ as above,
\[
S(r,\gamma,\ff_\ell) = \left\{
\begin{array}{ll}
\ell^{2r^2} & \gamma = 1 \\
\ell^{r^2-r} \frac{\#\sp_{2r}(\ff_\ell)}{\#\gl_r(\ff_\ell)} & \gamma
\not = 1
\end{array}
\right.
.
\]
\end{lemma}

\begin{proof}
 These computations use the following chain of standard
observations \cite{chav}.  Any characteristic polynomial of an element in
$\gsp_{2r}(\ff_\ell)$ is the characteristic polynomial of some semisimple
element $A$.  Moreover, the number of elements with such a
characteristic polynomial is
\[
\ell^{\dim Z(A) - r}\frac{\#\sp_{2r}(\ff_\ell)}{\#Z(A)(\ff_\ell)},
\]
where $Z(A)$ is the group of elements of $\sp_{2r}$ which commute
(inside $\gsp_{2r}$) with $A$.  From this, the computation of $S$
immediately follows.  Indeed, $A = \diag(1, \cdots, 1, \gamma, \cdots,
\gamma)$ is the unique semisimple element with characteristic
polynomial $(x-1)^r(x-\gamma)^r$.  If $\gamma = 1$, then $Z(A) =
Z(\id) = \sp_{2r}$, and $\dim Z(A) = 2r^2+r$.  If $\gamma\not = 1$,
then the centralizer $Z(A)$ is $\left\{ \twmt M00{(M\inv)^T} \right\} \iso
\gl_r$, a group of dimension $r^2$.  In either case, the lemma now
follows.
\end{proof}

\begin{lemma}\label{trec} With notation as above, $T(1,\gamma,\ff_\ell) =
S(1,\gamma,\ff_\ell)$.  For $g \ge 2$,
\begin{multline*}
T(g,\gamma,\ff_\ell) \\
= \sum_{\substack{r+s=g \\ 1 \le r \le g}}
\frac{\#\sp_{2g}(\ff_\ell)}{\#\sp_{2r}(\ff_\ell) \#\sp_{2s}(\ff_\ell)}
S(r,\gamma,\ff_\ell)
(\#\sp_{2s}(\ff_\ell) - T(s,\gamma,\ff_\ell)).
\end{multline*}
\end{lemma}

\begin{proof}
 The first claim is a tautology. For dimension $g\ge 2$, we
enumerate elements of $\gsp^\gamma_{2g}(\ff_\ell)$ which have one as
an eigenvalue.  First, we index elements $A$ of $W_{\eigenone, g}^\gamma(\ff_\ell)$ by
$r$, the order of vanishing of $f_A$ at $1$ if $\gamma \not = 1$, and
half that multiplicity if $\gamma = 1$.  To such an $A$ corresponds a
decomposition of $V$ as $U_{\eigenone, r}\oplus U_{\notone,s}$, where
$U_{\eigenone, r}$ and $U_{\notone, s}$ are
symplectic subspaces of dimensions $2r$ and $2s$, respectively;
$f_{A\rest{U_{\eigenone,r}}}(x)  = (x-1)^r(x-\gamma)^r$; and
$f_{A\rest{U_{\notone,s}}}(1)
\not = 0$.

The factor
$\frac{\#\sp_{2g}(\ff_\ell)}{\#\sp_{2r}(\ff_\ell)
\#\sp_{2s}(\ff_\ell)}$ counts the number of ways of decomposing $V =
U_{\eigenone,r}\oplus U_{\notone,s}$.  The penultimate factor
$S(r,\gamma,\ff_\ell)$ counts
the possibilities for $A$ acting on $U_r$, and the last factor
enumerates all choices for $A\rest {U_{\notone,s}}$.
\end{proof}

Roughly speaking, $\#W_{\eigenone,g}^\gamma/\#\sp_{2g}(\ff_\ell)$ is
about $\oneover \ell$.  In
fact, an argument similar to (but easier than) \ref{eigenweird} shows
that this ratio is between $(\ell/(\ell+1))^{2g^2+g}\oneover \ell$ and
$(\ell/(\ell-1))^{2g^2+g}\oneover\ell$.
Still, a more precise estimate isn't too difficult.

\begin{lemma}\label{eigenone} For each $g \ge 1$ there is a constant $c(g)$ such that
\[
\abs{\frac{T(g,\gamma,\ff_\ell)}{\#\sp_{2g}(\ff_\ell)} -
\tau^\gamma_{\eigenone,g}} \le c(g)(\tau^\gamma_{\eigenone,g})^3,
\]
where
\[
\tau^\gamma_{\eigenone,g} = \left\{
\begin{array}{ll}
\oneover{\ell-1} & \gamma\not = 1\\
\frac{\ell}{\ell^2-1} & \gamma = 1
\end{array}
\right..
\]
\end{lemma}

\begin{proof}
We treat the case $\gamma\not = 1$, and leave the remaining case for the
industrious reader.  Lemma \ref{trec} shows that
$T(1,\gamma,\ff_\ell)/\#\sp_2(\ff_\ell) = \tau^\gamma_1$, and that
$\abs{T(2,\gamma, \ff_\ell)/\#\sp_4(\ff_\ell)} =
(\ell^2-2)/((\ell-1)^2(\ell+1))$.  For $g \ge 3$,
\begin{eqnarray*}
\frac{T(g,\gamma,\ff_\ell)}{\#\sp_{2g}(\ff_\ell)} & = &
\sum_{\substack{r+s=g \\ 1 \le r \le g}}
\frac{S(r,\gamma,\ff_\ell)}{\#\sp_{2r}(\ff_\ell)}
\frac{\#\sp_{2s}(\ff_\ell) -
T(s,\gamma,\ff_\ell)}{\#\sp_{2s}(\ff_\ell)} \\
& = & \sum_{\substack{r+s = g\\ 1 \le r \le
g-1}}\frac{\ell^{r^2-r}}{\#\gl_r(\ff_\ell)} \left(1 -
\frac{T(s,\gamma,\ff_\ell)}{\#\sp_{2s}(\ff_\ell)}\right) +
\frac{\ell^{g^2-g}}{\#\gl_g(\ff_\ell)}.
\end{eqnarray*}
Now, for any $j\ge 1$, $\ell^j - 1 \ge (\ell-1) \ell^{j-1}$.  Thus,
\begin{eqnarray*}
\frac{\ell^{r^2-r}}{\#\gl_r(\ff_\ell)} & = &
\frac{ \ell^{r^2-r}}{\ell^{\half{r(r-1)}}\prod_{j=1}^r (\ell^j - 1)} \\
& \le & \frac{\ell^{\half{r(r-1)}}}{\prod_{j=1}^r \ell^{j-1}(\ell-1)}
\\
& = & \oneover{(\ell-1)^r}.
\end{eqnarray*}
Higher order terms -- those coming from $r > 2$ -- contribute
less than $O((\tau^\gamma_{\eigenone,g})^3)$ to
$T(g,\gamma,\ff_\ell)$; the lemma
is proved.
\end{proof}

\subsection{An intricate condition}\label{secweird}

Fix as before a dimension $g$, and consider $W_{\eigenweird,g}^\gamma\subset
\gsp_{2g}^\gamma(\ff_\ell)$, the
set of all elements $A$ whose characteristic polynomial $f_A(x)$
satisfies the following condition:

\begin{center}
\eigenweird \ \ \
\parbox{3in}{
Pairs of distinct roots of $f(x)$ over $\bar\ff_\ell$ do not multiply to $1$;
$f(x)$ has at most a single root  at $-1$; and $f(x)$ has at most a
double root at $1$.
}
\end{center}

While $W^\gamma_{\eigenweird,g}$ is presumably amenable to analysis in
the style of Lemma
\ref{eigenone}, we content ourselves with the following, coarser
estimate.

\begin{lemma}
\label{eigenweird} If $\gamma\not = 1$, then there is a constant $C(g)$ depending
only on $g$ such that
\[
\frac{\#W^\gamma_{\eigenweird,g}(\ff_\ell)}{\#\sp_{2g}(\ff_\ell)} \ge
\left(1 - \frac{C(g)}{\ell}\right)\left(\frac{\ell}{\ell+1}\right)^{2g^2+g}.
\]
\end{lemma}

\begin{proof}
Fix a $\gamma \in \ff_\ell\units$.  In fact, assume $\gamma\not = 1$;
for otherwise, $W_{\eigenweird,g}^\gamma = \emptyset$.  Consider the
space $\Xi^\gamma_g \iso \aff^g$ of characteristic polynomials of elements of
$\gsp_{2g}^\gamma$.
By considering successively the requirements for a point in
$\Xi^\gamma_g$ to satisfy $\eigenweird$, we will show that
$\eigenweird$ is a Zariski open condition.

The first condition is that $f(x)$ and $f(1/x)$ have no common
root.  This is clearly an open condition, as it is equivalent to the
disjointness of $\spec \frac{\ff_\ell[x]}{f(x)}$ and $\spec
\frac{\ff_\ell[x]}{f(1/x)}$ inside $\spec \ff_\ell[x] \iso
\aff^1_{\ff_\ell}$.  The second condition says that at least one of $f(-1)$
and $f'(-1)$ is nonzero; and the final condition says that at least
one of $f(1)$, $f'(1)$, and $f''(1)$ is nonzero.

It is clear that $\Psi^\gamma_{\eigenweird,g}$ is nonempty if and only
if $\gamma\not = 1$.  So there is a constant $C(g)$, depending on $g$
but not on $\ell$, such that
\[
\#\Psi^\gamma_{\eigenweird,g}(\ff_\ell) \ge \ell^g - C(g)\ell^{g-1}.
\]
Invoking \ref{psitow} now proves the lemma.
\end{proof}

\section{\fm\  holds generically}

Following the abstract situation at the end of Section \ref{mono}, say
that an abelian variety $X$ over a finite field $k$ has $\notone$ if
it {\em does not} have a rational $\ell$-torsion point over $k$.  Say
that $X/k$ has $\eigenweird$ if its characteristic polynomial of
Frobenius, taken modulo $\ell$, satisfies $\eigenweird$ as in
Section \ref{secweird}.
Recall that a curve has $\notone$ or $\eigenweird$ if its Jacobian
does, and that a function field has such a property if its proper,
smooth model does.

\begin{theorem}\label{oldmainth} If a function field satisfies $\notone$ or
$\eigenweird$, then it satisfies \fm.
\end{theorem}

\begin{proof} If the function field satisfies \notone, then we know that
$\ell$ does not divide the class number of the function field.  We
thus have the conditions of Theorem~2 of~\cite{Holden99}, which shows
that the function field satisfies \fm.
On the other hand if the function field satisfies $\eigenweird$, then
we are in the situation of Remark~4.10 of~\cite{Holden99}, and again
the function field satisfies \fm.
\end{proof}

We will now see that, in some sense, most function fields fall under
the aegis of \ref{oldmainth}.

\begin{theorem}\label{mainth} Let $\calc\ra\calm \ra \spec k_0$ be a smooth,
irreducible, proper family of proper curves of genus $g \ge 1$ with
full $\ell$-monodromy.  Fix a
$\gamma\in \ff_\ell\units$.  Let $\st{k_n}$ be a collection of
extensions of $k_0$ such that $\lim_{n\ra\infty} \#k_n = \infty$, and,
for all $n$, $\gamma(k_n) = \gamma$. For $\bullet =
\notone,\eigenweird$,
\[
\lim_{n\ra\infty}
\frac{\#\st{x\in \calm(k_n) : \calp^\gamma_\bullet(\calc_x)}}
{\#\calm(k_n)}
=
\frac{\#W^\gamma_{\bullet,g}(\ff_\ell)}{\#\sp_{2g}(\ff_\ell)}
=
\sigma^\gamma_{\bullet,g}
.
\]
There are constants $C^\gamma_{\bullet,g}$, independent of $\ell$,
such that:
\begin{eqnarray*}
\sigma^1_{\eigenweird,g} & = & 0; \\
\sigma^\gamma_{\eigenweird,g} &\ge&\left( 1 -
\frac{C^\gamma_{\eigenweird,g}}{\ell}\right)\left(
  \frac{\ell}{\ell+1}\right)^{2g^2+g}\text{ if }\gamma\not = 1;\\
\abs{ \sigma^\gamma_{\notone,g} - \tau^\gamma_{\notone,g}} & \le &
C^\gamma_{\notone,g}(1-\tau^\gamma_{\notone,g})^3 =
O\left(\oneover{\ell^3}\right);\\
\tau^1_{\notone, g} & = & 1 - \frac{\ell}{\ell^2-1}; \\
\tau^\gamma_{\notone,g} & = & 1 - \oneover{\ell-1} \text{ if
  }\gamma\not = 1.
\end{eqnarray*}
\end{theorem}

\begin{proof} By \ref{katz}, the proportion of curves with either property
$\notone$ or $\eigenweird$ converges to $\sigma^\gamma_{\bullet,g}$,
the appropriate proportion of elements of
$\gsp^\gamma_{2g}(\ff_\ell)$.  Lemmas \ref{eigenone} and
\ref{eigenweird} estimate these values for $\notone$ and
$\eigenweird$, respectively.
\end{proof}

The exact same techniques let us compute the proportion of abelian
varieties over a finite field which have a rational $\ell$-torsion
point.  For a natural number $N$, let $\cala_{g,N}$ denote the fine
moduli scheme of triples $(A,\lambda,\phi)$ consisting of an abelian
scheme $A$, a principal polarization $\lambda$, and symplectic
principal level $N$
structure $\phi$.

\begin{prop}\label{abvar} Let $N = \ell N' \ge 3$ be a natural number relatively prime to
$p$, and let $k_0$ be a finite field of characteristic $p$
containing a primitive $N^{th}$ root of unity.  Let $\st{k_n}$ be a
tower of extensions of $k_0$ such $\lim_{n\ra\infty}\#k_n = \infty$
and, for $n$ sufficiently large, $\#k_n \equiv 1 \bmod \ell$.  Then
\[
\lim_{n\ra\infty}
\frac{\#\st{(A,\lambda,\phi)\in \cala_{g,N}(k_n): A[\ell](k_n)\supsetneq
\st 1 }}{\#A_{g,N}(k_n)}
=
\frac{\ell}{\ell^2-1} + O(\oneover{\ell^3}),
\]
where the constant in the error term $O(\oneover{\ell^3})$ depends
only on $g$.
\end{prop}

\begin{proof}
The proof of \ref{mainth} uses only statements about abelian
varieties, and thus applies in this setting, too.  Let
$\calg/\cala_{g,N}$ be the sheaf of $\ell$-torsion of the universal
abelian variety over $\cala_{g,N}$.  Fix a geometric point $\bar x \in
\cala_{g,N}$ which is the Jacobian of a curve.  Since the Torelli
locus already has full monodromy (\ref{full}), the image of
$\pi_1^\geom(\cala_{g,N}, \bar x )
\ra \calg_{\bar x}$ is $\sp_{2g}(\ff_\ell)$.  Thus, all the machinery
exposed in this paper applies, and the result follows.  By \ref{katz},
as $\#k_n \ra \infty$ the proportion of abelian varieties with a
rational $\ell$-torsion point approaches the proportion of symplectic
matrices with one as an eigenvalue.  The latter ratio, or at least its
leading term, is computed in \ref{eigenone}.
\end{proof}

\begin{remark}
  We would like to comment briefly on the collections of fields
$\st{k_n}$ in \ref{mainth} and \ref{abvar}.  On one hand, the
collection of extensions $k$ of $k_0$ with $\#k \equiv 1 \bmod \ell$
is cofinal, in that any extension of $k_0$ is a subfield of such a
subfield.  Thus, it seems natural to take limits over towers of such
fields; this explains our choice in \ref{abvar}.  (It is not hard to
adapt the statement for a different $\gamma(k_n)$.)

On the other hand, if $\gamma,\gamma'\not=1$, then
$W^\gamma_{\notone,g} = W^{\gamma'}_{\notone,g}$.  Thus, in
\ref{mainth}, the collection of fields $\st{k_n}$ may
be generalized to any collection of increasingly large finite
extensions of $k_0$, so long as $\gamma(k_n)\not = 1$ for sufficiently
large $\gamma$.
\end{remark}

\begin{remark} Theorem \ref{mainth} shows that approximately
$\frac{\ell-2}{\ell-1}$ of all function fields satisfy $\notone$, and thus
\fm.  Similarly, roughly $\frac{\ell-1}{\ell}$ of all function fields
satisfy $\eigenweird$.  In fact, it seems likely that on the order of $1 -
\oneover{\ell^2}$ function fields satisfy $\notone$ or $\eigenweird$, and
that at least this proportion of function fields has \fm.  (One could
compute this number directly, but at present the relatively modest payoff
does not seem to warrant the detailed combinatorics required.)  To see
this, we argue on the level of characteristic polynomials, using Lemma
\ref{psitow} to help us pass from characteristic polynomials to group
elements.  In \ref{eigenweird} we showed that $\eigenweird$ is a Zariski
open condition on $\Xi^\gamma_g$.  Similarly, one directly sees that
$\eigenone$ is a closed condition on $\Xi^\gamma_g$ since $f(x) \in
\Psi^\gamma_{\eigenone,g}$ if and only if $f(1) = 0$ and that $\notone$ is
Zariski open.  The closed conditions which trace the complement of
$\Psi_{\eigenweird,g}^\gamma$ and $\Psi_{\notone,g}^\gamma$ are
independent, and $\Psi^\gamma_{\notone \text{ or }\eigenweird, g}$ is the
complement of a closed set of codimension two in $\Xi^\gamma_g$.  If
proportions of characteristic polynomials are directly reflected in
proportions of elements of the symplectic group, then about
$1-\oneover{\ell^2}$ function fields satisfy $\notone$ or $\eigenweird$.
\end{remark}

% \bibliographystyle{abbrv}
% \bibliography{josh}

% the following is from the fm-jtnbv2.bbl file

\end{document}